\documentclass[twoside,12pt]{article}
\usepackage{amsfonts}
\usepackage{amsmath}
\usepackage{amssymb}
\usepackage[mathscr]{eucal}
\input epsf

\newtheorem{THEOREM}{Theorem}
\newtheorem{LEMMA}[THEOREM]{Lemma}
\newtheorem{DEFINITION}[THEOREM]{Definition}
\newtheorem{REMARK}[THEOREM]{Remark}
\newtheorem{PROPOSITION}[THEOREM]{Proposition}
\newtheorem{COROLLARY}[THEOREM]{Corollary}

\def\beginrem{\begin{REMARK}\rm}
\def\endrem{\end{REMARK}}
\def\beginthm{\begin{THEOREM}}
\def\endthm{\end{THEOREM}}
\def\beginlem{\begin{LEMMA}}
\def\endlem{\end{LEMMA}}
\def\beginprop{\begin{PROPOSITION}}
\def\endprop{\end{PROPOSITION}}
\def\begincor{\begin{COROLLARY}}
\def\endcor{\end{COROLLARY}}
\def\begindef{\begin{DEFINITION}\rm}
\def\enddef{\end{DEFINITION}}

\newcommand{\cab}{,\allowbreak}

\makeatletter
\newcommand{\gclass}{\relax\ifmmode\mathbb G\else
$\m@th\mathbb G$\fi}
\newcommand{\fclass}{\relax\ifmmode\mathbb F\else
$\m@th\mathbb F$\fi}

\newcommand{\nclass}{\relax\ifmmode \mathbb N\else
$\m@th\mathbb N$\fi}

\newcommand{\abs}[1]{\relax\ifmmode{|}#1{|}\else
$\m@th{|}#1{|}$\fi}
\newcommand{\cardinal}[1]{\relax\ifmmode\left|#1\right|\else
$\m@th\left|#1\right|$\fi}

\newcommand{\smalldisplayspace}
{\abovedisplayskip = 1.5ex plus .15ex minus .3ex
\belowdisplayskip = 1.5ex plus .15ex minus .3ex }

\newcommand{\bigspace}{\par\ifdim\lastskip<2 ex 
\removelastskip\penalty-150 \vskip 2 ex plus .25ex minus .25ex\fi}

\newcommand{\paraspace}{\par\ifdim\lastskip<1.5 ex \removelastskip\penalty-100 
\vskip 1.5 ex plus .25ex minus .25ex\fi}

\newcommand{\smallspace}{\par\ifdim\lastskip<1 ex \removelastskip\penalty-100 
\vskip 1 ex plus .2ex minus .2ex\fi}

\newcommand{\proofspace}{\par\ifdim\lastskip<2 ex
\removelastskip\penalty-100  \vskip 2 ex plus .25ex minus .25ex\fi}

\def\endproofsymbol{\mbox{$\square$}}

\def\endproof{\Tag{\endproofsymbol}\proofspace}
\newcommand{\proof}{\par\noindent{\bf Proof.}\kern.5em}

\newcommand{\Tag}[1]{\ifvmode\else\unskip\fi
\nobreak\hfil\penalty50 \hskip2em \null
\nobreak\hfil#1\skip@\parfillskip\parfillskip\z@skip
\count@\finalhyphendemerits\finalhyphendemerits\z@\par
\parfillskip\skip@\finalhyphendemerits\count@}

\newdimen\LabeLmargin\LabeLmargin=0pt
\newdimen\LabeLwidth\LabeLwidth=0pt
\newdimen\Hangamount\Hangamount=0pt
\newdimen \ListSpace \ListSpace=2ex
\newdimen\LabelSpace \LabelSpace=1ex

\def\setListSpace{\ifdim\lastskip<\ListSpace\removelastskip\penalty-100
\vskip \ListSpace plus .15 \ListSpace minus .15 \ListSpace\fi} 

\def\BeginList#1#2{\par\begingroup\setListSpace
\parindent=0pt \LabeLmargin=\Hangamount\setbox0=\hbox{#1}
\advance\LabeLmargin by\wd0 \Hangamount=\LabeLmargin
\setbox0=\hbox{#2}\advance\Hangamount by\wd0 \LabeLwidth=\wd0}

\def\setLabelSpace{\ifdim\lastskip<\ListSpace
\removelastskip\penalty-75  \vskip \LabelSpace plus .15 \LabelSpace
minus .15 \LabelSpace\fi}

\def\LabeL#1{\par\setLabelSpace \hangindent=\Hangamount
\hangafter=1 \hskip\LabeLmargin\hbox to\LabeLwidth{#1\hfill}\ignorespaces}

\def\@begintheorem#1#2{\trivlist\item[\hskip\labelsep{\bfseries
#1\ #2.}]\itshape}

\makeatother

\begin{document}\thispagestyle{empty}

\begin{center}\large\bf A Note on Ternary Sequences of Strings of 0 and 1
\end{center}
  {\tabskip 0pt plus 1fil \halign to \hsize{#\hfil\cr
A. R. Mehta\cr
Department of Mathematics\cr
Indian Institute of Technology\cr
Guwahati 781\,039\cr
\noalign{\vskip 1ex}
Email: mehta@iitg.ernet.in\cr
\noalign{\vskip 2ex}
G. R. Vijayakumar\cr
School of Mathematics\cr
 Tata Institute of Fundamental Research\cr
 Homi Bhabha Road, Colaba\cr
 Mumbai 400\,005\cr
	\noalign{\vskip 1ex}
 Email: vijay@math.tifr.res.in\cr}}

 \begin{abstract}\noindent 
B. D. Acharya has conjectured  that if
$\bigl(A_i: i=1\cab 2\cab \ldots \cab 2^{\abs X}-1\bigr)$ is a
permutation of all nonempty subsets of a set $X$ with at least two
elements such that for each even positive integer $j<2^{\abs X}-1$,
$A_{j-1}\triangle A_j\triangle A_{j+1}=\emptyset$, then $\abs X=2$.
In this article, we show that if the cardinality of a set $X$ is
 more than four, then  a permutation as described above indeed
exists.
\end{abstract}
 \smallspace\noindent
 \paraspace
  
\noindent Throughout this article, $\gclass$ is the group defined on
$\{0,1\}$. Let $n$ be any 
positive integer. We denote the identity of $\gclass^n$ usually by 0 and
sometimes by $z_n$ to avoid ambiguity.
 Let $\fclass^n$ be the
set of all nonzero elements in $\gclass^n$; when $n\neq 1$,
a permutation $\left(v_i: i=1\cab 2\cab \ldots,
2^n-1\right)$ of $\fclass^n$ is called {\it ternary\/} if 
for each $i\in \{2\cab 4\cab\ldots \cab 2^n-2\}$,
$v_{i-1}+v_i+v_{i+1}=0$; when such a permutation
exists,  $\gclass^n$ is called {\it sequentially
ternary}. In this article, elements in $\gclass^n$ are represented
by  strings of `symbols' where each symbol is an element in a power
of $\gclass$. For example, if $\alpha=01$ and $\beta=10$, then
$1\alpha0\beta$ is the element $(1\cab 0\cab 1\cab 0\cab 1\cab 0)$ in
$\gclass^6$. 
\paraspace
In \cite{bda}, B. D. Acharya has observed that $\gclass^3$ is not
sequentially ternary and  conjectured that if $n$ is an
integer which is larger than 2, then $\gclass^n$ is
{\it not\/} sequentially ternary. This article is an outcome of 
settling this conjecture. For basic group theoretic results needed in
this connection, we rely on \cite{her}.

\begin{table}\vbox to \vsize{\epsfysize=.62\vsize
\centerline{\epsffile[43 38 794 830]{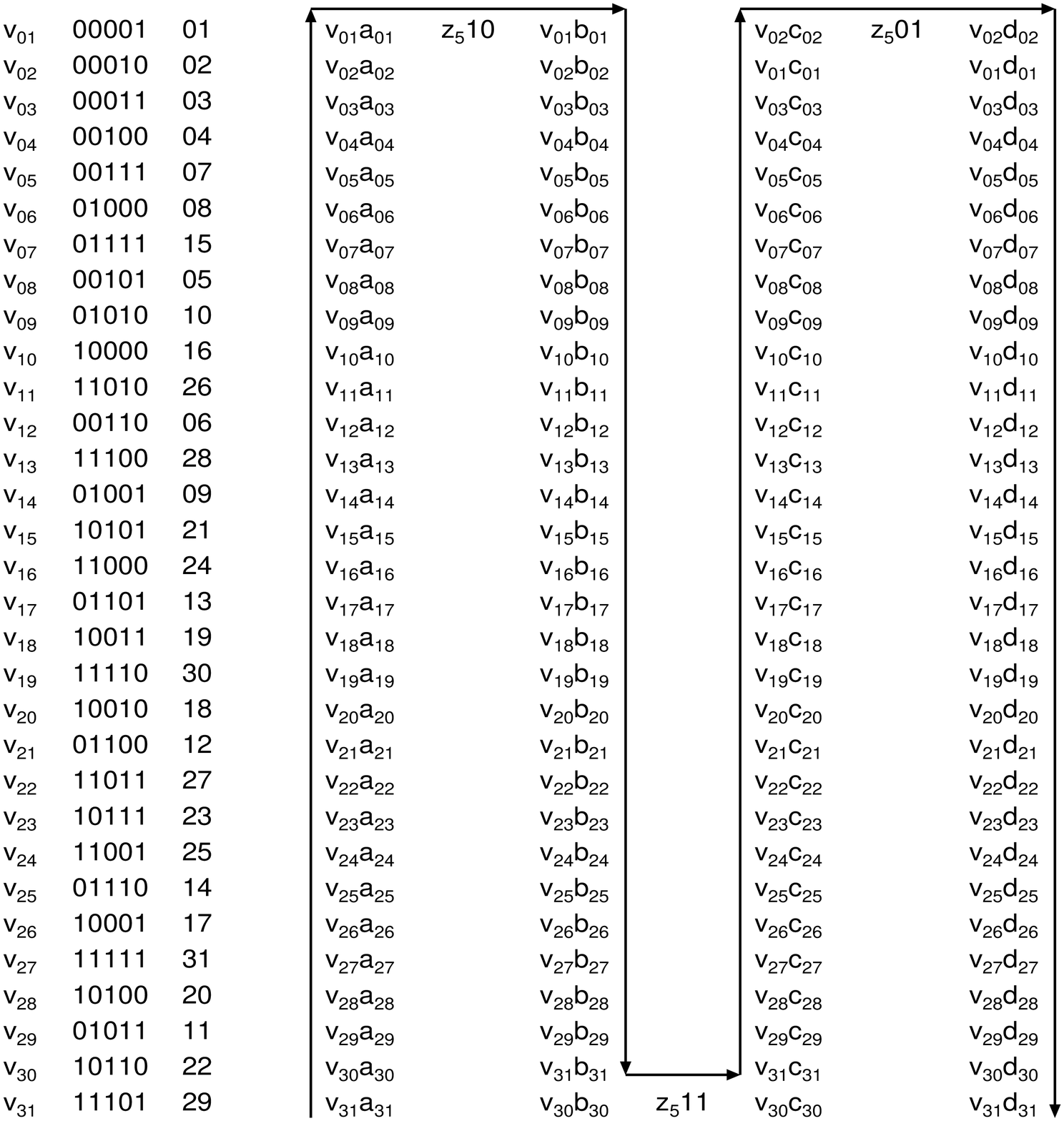}}
\bigspace \epsfysize=.31\vsize
 \centerline{\epsffile[68 430 830 825]{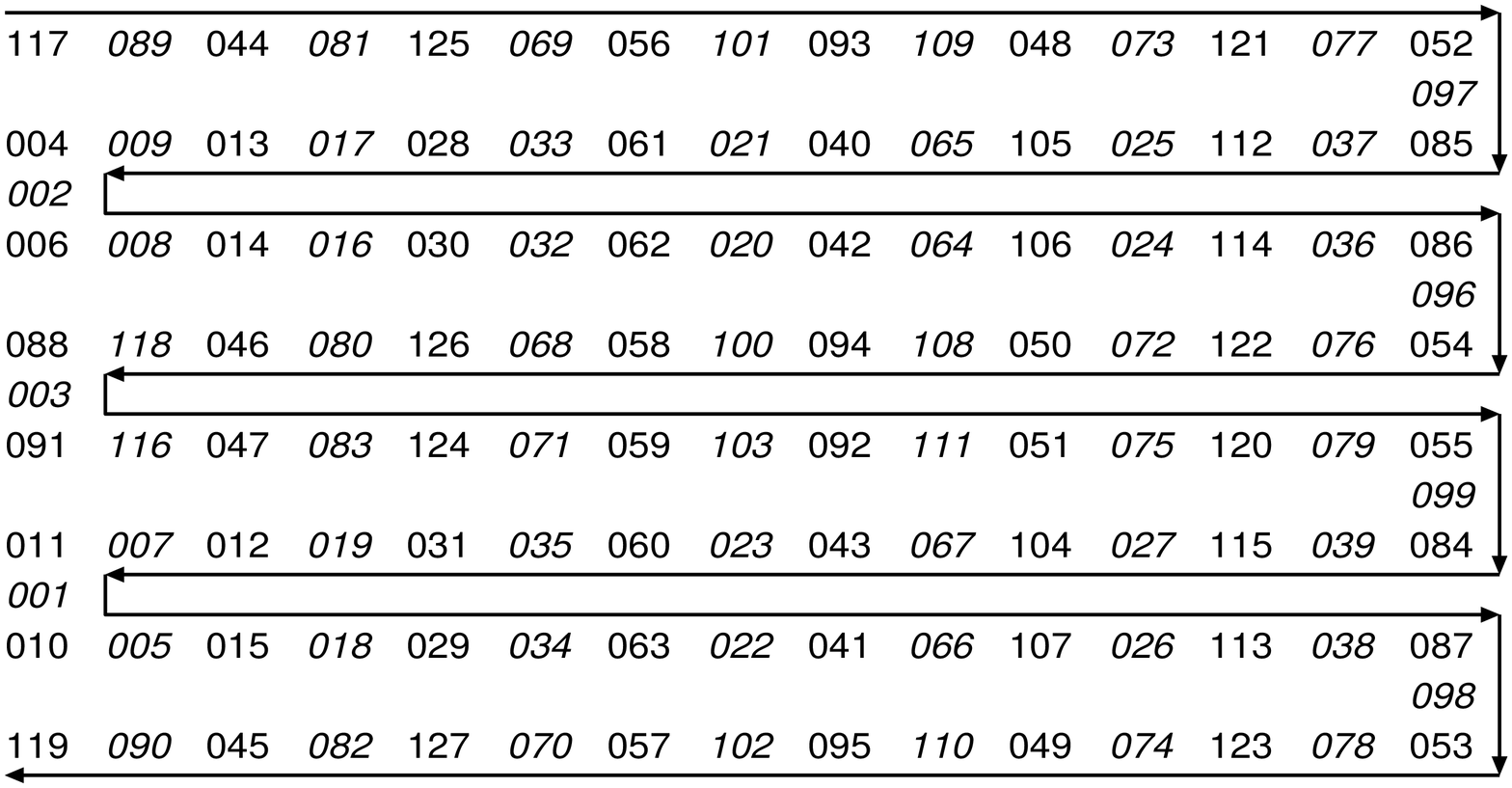}}
\paraspace \centerline{Figure 1}\vfil}
\end{table}

\beginprop\label{recursive} For any 
integer $n>2$, if $\gclass^n$ is sequentially ternary,
then $\gclass^{n+2}$ is also sequentially ternary.
\endprop
\proof Let  $k=2^n-1$. Define sequences
$\bigl(a_i\bigr)_{i=1}^k$, $\bigl(b_i\bigr)_{i=1}^k$,
$\bigl(c_i\bigr)_{i=1}^k$ and  $\bigl(d_i\bigr)_{i=1}^k$ in $\gclass^2$
as described below.
\paraspace {\tabskip 0pt plus 1fil

\newcommand\tempcs[1]{$\left\lbrace\vcenter{\vbox{\ialign{##\quad
&##\hfil\cr#1\crcr}}}\right.$}

\halign to \hsize{\hfil$#={}$\tabskip 0pt&#\tabskip 0pt plus 1fil\cr
a_i&\tempcs{00,&if $i\equiv1$ (mod 4);\cr
01,&otherwise.\cr}\cr\noalign{\vskip 1ex}
b_i&\tempcs{10,&if $i\equiv1$ (mod 2);\cr
00,&otherwise.\cr}\cr\noalign{\vskip 1ex}
c_i&\tempcs{00,&if $i\equiv3$ (mod 4);\cr
11,&otherwise.\cr}\cr\noalign{\vskip 1ex}
d_i&\tempcs{01,&if $i\equiv1$ (mod 4);\cr
11,&if $i\equiv3$ (mod 4);\cr
10,&otherwise.\cr}\cr}
\bigspace}
\noindent It is easy to see that\smallspace
\noindent (1) for each $i\in \{1\cab 2\cab \ldots \cab k\}$ the
strings $a_i$, $b_i$, $c_i$ and $d_i$ are  distinct and\vskip 1ex
\noindent (2) for each $i\in \{2\cab 4\cab \ldots \cab k-1\}$,
$$\abovedisplayskip 1ex \belowdisplayskip 0pt
a_{i-1}+a_i+a_{i+1}=b_{i-1}+b_i+b_{i+1}=c_{i-1}+c_i+c_{i+1}=
d_{i-1}+d_i+d_{i+1}=z_2.$$
\paraspace

\noindent Let
 $\bigl(v_i: i=1\cab\ldots\cab k\bigr)$  be a ternary
permutation of $\fclass^n$.
 Define a sequence $\bigl(w_i\bigr)_{i=1}^{4k+3}$ in
 $\fclass^{n+2}$ as described below. \paraspace
 {\leftskip =3em \parindent =0pt \rightskip 0pt plus 1fil
For each $i\in \{1\cab 2\cab \ldots, k\}$,
$w_i=v_{k-i+1}a_{k-i+1}$
and $w_{k+1}=z_n10$. \smallspace
 For each $i\in \{k+2\cab k+3\cab \ldots \cab
2k-1\}$, $w_i=v_{i-k-1}b_{i-k-1}$. \smallspace
 $w_{2k}=v_kb_k$, $w_{2k+1}=v_{k-1}b_{k-1}$,
$w_{2k+2}=z_n11$, $w_{2k+3}=v_{k-1}c_{k-1}$ and
$w_{2k+4}=v_kc_k$.
\smallspace
For each $i\in \{2k+5\cab 2k+6\cab \ldots \cab 3k\}$,
$w_i=v_{3k+3-i}c_{3k+3-i}$.\smallspace 
$w_{3k+1}=v_1c_1,\ w_{3k+2}=v_2c_2$,
$w_{3k+3}=z_n01$, $w_{3k+4}=v_2d_2$ and
$w_{3k+5}=v_1d_1$.
\smallspace \noindent
For each $i\in\{3k+6\cab 3k+7\cab \ldots \cab 4k+3\}$,
$w_i=v_{i-3k-3}d_{i-3k-3}$.\paraspace}
\noindent 
It can be verified that
$\{w_i: i=1\cab 2\cab \ldots \cab
4k+3\}=\{v_ia_i\cab v_ib_i\cab v_ic_i\cab v_id_i:i=1\cab \ldots \cab k\}
\cup\{z_n10\cab z_n11\cab z_n01\}$. Therefore by (1), all terms of 
$\bigl(w_i\bigr)_{i=1}^{4k+3}$ are distinct. 
From the definition of this sequence we have the following.
\paraspace {\tabskip 0pt plus 1fil
\halign to\hsize{\hfil$#$\tabskip 0pt&${}=#$\hfil
\tabskip 0pt plus 1fil\cr
w_k+w_{k+1}+w_{k+2}&v_1a_1+z_n10+v_1b_1;\cr
w_{2k+1}+w_{2k+2}+w_{2k+3}&v_{k-1}b_{k-1}+z_n11+v_{k-1}c_{k-1};\cr
\mathrm{and}\ w_{3k+2}+w_{3k+3}+w_{3k+4}&v_2c_2+z_n01+v_2d_2.\cr}}
\paraspace 
\noindent It is easy to see that  for each
of the above three equations, its right side is zero;
 therefore for all $i\in\{k+1\cab
2k+2\cab 3k+3\}$, $w_{i-1}+w_i+w_{i+1}=0$; by using (2) it can be
easily verified that for each  $i\in \{2\cab 4\cab 6\cab \ldots
\cab 4k+2\}\setminus \{k+1\cab 2k+2\cab 3k+3\}$
 also, the just mentioned equality holds. Therefore,
the permutation $\bigl(w_i: i=1\cab \ldots\cab 4k+3\bigr)$ is ternary.
\endproof 

\begin{table}\vbox{\epsfxsize=.95\hsize
 \centerline{\epsffile[39 221 791 563]{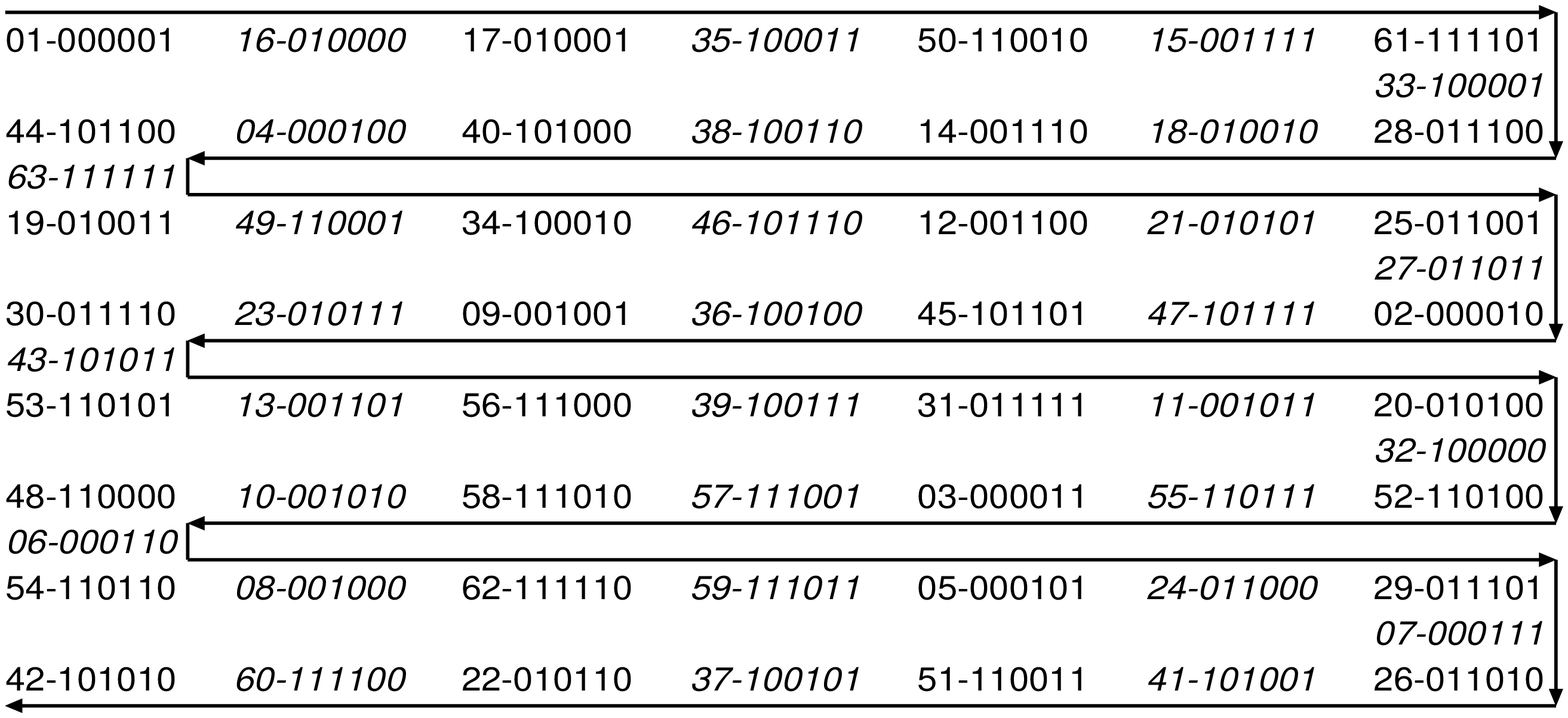}}}
 \paraspace \centerline{Figure 2}
\end{table}

\beginrem\label{string-5-7}  Figure 1
 illustrates the method described by Proposition \ref{recursive}
for $n=5$. It can be verified
that the sequence of strings of length 5 which
is displayed in this figure---for each $k\in \{1\cab \ldots
\cab 31\}$, $k$-th string is labelled by $v_{0k}$ or $v_k$ and
treating each string as 
a binary number, its value in the decimal form is written on
its right---is a ternary permutation of $\fclass^5$.
The  ordering of the other two sequences in the figure
are  indicated by  directed lines.
(Note that $v_{31}b_{31}$, $v_{31}c_{31}$, $v_{01}c_{01}$ and 
$v_{01}d_{01}$ appear in between $v_{29}b_{29}$ and $v_{30}b_{30}$,
$v_{30}c_{30}$ and $v_{29}c_{29}$, $v_{03}c_{03}$ and $v_{02}c_{02}$,
and $v_{02}d_{02}$ and $v_{03}d_{03}$ respectively.)
The ternary permutation of $\fclass^7$ thus obtained is represented by
 a sequence of decimal numbers. (The even numbered terms of the sequence
 are italicized.) 
\endrem
Let $n$ be an integer which exceeds 1.
 Let $u$ be any element in $\fclass^n$. It is easy to see that
$\gclass^n$ can be partitioned into $2^{n-1}$ pairs such that
sum of the elements of each pair is $u$; therefore the sum of all
elements in $\gclass^n$ is $0$. We use this (well known) fact to
settle the first part of the next result.
 \beginprop\label{string-3-4}
  Neither $\gclass^3$ nor\/ $\gclass^4$ is sequentially ternary.
 \endprop
\proof Let $\bigl(a_i:i=1\cab \ldots\cab 7\bigr)$ be a
permutation of $\fclass^3$. Since $a_1+a_2+\cdots +a_7=0$ 
either $a_1+a_2+a_3$ or $a_5+a_6+a_7$ is nonzero; therefore
this permutation is not ternary.
\paraspace Next, suppose that $\bigl(a_i:i=1\cab \ldots \cab
15\bigr)$ is a ternary permutation of $\fclass^4$. Let $H$ be
the subgroup of $\gclass^4$
which is  generated by $\{a_2\cab a_3\cab a_4\}$. Obviously,
$a_1\cab a_5\in H$. Let us find the other two nonzero elements
in $H$. Noting that if $K$ is a subgroup of cardinality 4,
then $\abs{H\cap K}\geqslant 2$, because
$$\smalldisplayspace
\abs{H\cap K}=\frac{\abs H\times \abs K}{\abs{HK}}\geqslant
\frac{\abs H\times 4}{\abs{\gclass^4}}\geqslant 2,$$
 for each $i\in \{7\cab 9\cab 11\cab 13\}$, the fact that
 $\{0\cab a_i\cab a_{i+1}\cab a_{i+2}\}$ is a subgroup implies that
$H\cap\{a_i\cab a_{i+1}\cab a_{i+2}\}\ne \emptyset$; from this, we
find that $a_9\cab a_{13}\in H$. 
 Thus $\{0\cab a_1\cab a_2\cab a_3\cab a_4\cab a_5\cab a_9\cab
a_{13}\}$ is  a subgroup; by
symmetry, $\{0\cab a_{15}\cab a_{14}\cab a_{13}\cab 
a_{12}\cab a_{11}\cab a_7\cab a_{3}\}$ is also a subgroup.
The cardinality of the intersection of these two subgroups
is 3---a contradiction.
\endproof

\beginrem\label{string-6} In Figure 2, a sequence of 63 terms is displayed;
each term has an element of $\fclass^6$ and its value in decimal form.
It can be verified that the  strings of this sequence
form a ternary permutation of $\fclass^6$.
\endrem
 Now combining the
information we have from Propositions \ref{recursive}
and \ref{string-3-4} and Remarks \ref{string-5-7} and \ref{string-6}, we get
the following.\paraspace

\noindent {\bf Theorem.}\enspace{\it For any integer $n\geqslant
2$, $\gclass^n$ is sequentially ternary if and only if $n$ is neither
3 nor 4.}


\begin{thebibliography}{9}
 
\bibitem{bda} B. D. Acharya, Set valuations of graphs and their
applications, {\it MRI Lecture Notes in Applied Mathematics}, Vol. 2,
 Mehta Research Institute, Allahabad (1983).
 \bibitem{her} I. N. Herstein, {\it Abstract Algebra}, $3^{\mathrm{rd}}$
edition, Prentice Hall, Upper Saddle River, New Jersey (1996).
 \end{thebibliography}
\end{document}